\documentclass[a4paper]{amsart}
\title[Derived Schwarz map]
{%
  Derived Schwarz map\\
  of the hypergeometric differential equation \\
  and a parallel family of flat fronts
}
\date{February 28, 2007}
\usepackage[dvips]{graphicx} 
\usepackage{amssymb}
\usepackage{rsfs}
\usepackage[usenames]{color}

\usepackage{amsthm}
\theoremstyle{plain}
 \newtheorem{theorem}{Theorem}[section]
 \newtheorem{proposition}[theorem]{Proposition}
 \newtheorem{lemma}[theorem]{Lemma}
 
\theoremstyle{remark}
 \newtheorem{remark}[theorem]{Remark}
 
\numberwithin{equation}{section}

\newcommand{\R}{\boldsymbol{R}}
\newcommand{\C}{\boldsymbol{C}}
\newcommand{\PP}{\boldsymbol{P}}
\newcommand{\BB}{\boldsymbol{B}}

\newcommand{\HH}{\boldsymbol{H}}

\def\transpose#1{\mathord{\mathopen{{\vphantom{#1}}^t}#1}} 
\def\qed{$\Box$\bigbreak} 
\renewcommand{\dfrac}[2]{%
  \frac{\displaystyle{#1}}{\displaystyle{#2}}}
\author{Takeshi Sasaki}
\address[Sasaki]{%
   Department of Mathematics,
   Kobe University,
   Kobe 657-8501, Japan%
}
\email{sasaki@math.kobe-u.ac.jp\qquad {\sl fax number:{\rm +81-78-803-5610}}}

\author{Kotaro Yamada}
\address[Yamada]{%
   Faculty of Mathematics,
   Kyushu University,
   Fukuoka 812-8581, Japan%
}
\email{kotaro@math.kyushu-u.ac.jp}

\author{Masaaki Yoshida}
\address[Yoshida]{%
   Faculty of Mathematics,
   Kyushu University,
   Fukuoka 810-8560, Japan%
}
\email{myoshida@math.kyushu-u.ac.jp}



\begin{document}
\maketitle

\begin{abstract}
In the paper \cite{SYY} we defined a map, called the hyperbolic
Schwarz map, from the one-dimensional projective space
to the three-dimensional hyperbolic space
by use of solutions of the hypergeometric differential equation, and thus
obtained closed flat surfaces belonging to the class of flat fronts.
We continue the study of such flat fronts in this paper.
First, we introduce the notion of derived Schwarz maps of the
hypergeometric differential equation and, second, we construct a
parallel family of flat fronts connecting the classical Schwarz map and
the derived Schwarz map.
\end{abstract}

\section{Introduction} 
Consider the {\it hypergeometric differential equation} 
\begin{equation}\tag*{$E(a,b,c)$}\label{eq:Eabc}
    x(1-x)u''+\{c-(a+b+1)x\}u'-abu=0,
\end{equation}
and define its {\it Schwarz map\/} 
as a multi-valued map on $X=\C-\{0,1\}$ by
\begin{equation}\tag{S}\label{eq:S}
    S:X \ni x\longmapsto u_0(x):u_1(x)\in \PP^1,
\end{equation}
where $u_0$ and $u_1$ are linearly independent solutions of
\ref{eq:Eabc}
and $\PP^1$ is the complex projective line.
A change of the unknown $u$ by multiplying a non-zero function, 
takes the equation into the SL-form: 
\begin{equation}\tag{SL}\label{eq:SL-form}
   u''-q(x) u = 0.
\end{equation}
The coefficient $q$ is expressed as
\[ q=-\{S;x\} \]
where $\{S;x\}$ is the Schwarzian derivative
 \[
    \{S;x\}:=\frac{1}{2}\left(\frac{S''}{S'}\right)'
             -\frac{1}{4}\left(\frac{S''}{S'}\right)^2.
 \]
For two linearly independent solutions $u_0$ and $u_1$ to this
equation, we define the {\it derived Schwarz map}  
\begin{equation}\tag{DS}\label{eq:DS}
   DS:X\ni x\longmapsto u'_0(x):u'_1(x)\in \PP^1,
\end{equation}
and the {\it hyperbolic Schwarz map}
\begin{equation}\tag{HS}\label{eq:HS}
   HS:X\ni x\longmapsto H(x)=U(x)\ {}^t\overline U(x)\in \HH^3,
\end{equation}
where 
\begin{equation}\label{eq:hol-lift}
    U=
      \left(\begin{array}{cc}u_0&u_0'\\u_1&u_1'\end{array}\right),
\end{equation}
and $\HH^3$ is the hyperbolic 3-space identified with the space of
positive $2\times2$-hermitian matrices modulo diagonal ones. 
The hyperbolic Schwarz map is considered as a flat front in $\HH^3$
in the sense of \cite{KUY}, that is, a flat surface of certain kind
of singularities.
See Section~\ref{sub:geom-1} and  Section~\ref{sub:flat-fronts}
for details.
We regard $\PP^1$ as the ideal boundary  $\partial\HH^3$ of $\HH^3$. 
Then $S$ and $DS\colon{}X\to \PP^1=\partial\HH^3$ are
considered as the two {\em hyperbolic Gauss maps\/} of the flat front $HS$.

We assume the parameters $a,b$ and $c$ are real and satisfy the condition
\begin{equation}\label{mule1}
   |\mu_0|,\quad|\mu_1|,\quad|\mu_\infty|<1,
\end{equation}
where 
\[
  \mu_0=1-c,\quad\mu_1=c-a-b,\quad\mu_\infty=b-a
\]
are exponent-differences at $0,1$ and $\infty$, respectively.
The Schwarz map gives a conformal equivalence between the upper half
part 
\[
  X_+:=\{x\in X=\C\setminus\{0,1\}\mid \Im(x)\ge0\}
\]
of $X$ and the image $T:=S(X_+)$, which is bounded by three arcs, 
a Schwarz triangle.  
Though the image $DS(X_+)$ is bounded by the three circles generated by
the three arcs bounding the Schwarz triangle, 
the situation depends on the parameters $(a,b,c)$. 
The image surface $HS(X_+)$ has, in general, singularities; 
the situation also depends on the parameters $(a,b,c)$. 
We study such dependence. 

On the other hand, 
there is a 1-parameter (parallel) family of surfaces 
(maps from $X$ to $\HH^3$) in the hyperbolic 3-space,
such that 
the Schwarz and the derived Schwarz maps are the two extremes 
(which have the images in $\PP^1$), and the hyperbolic Schwarz map 
is a generic member. 
For a typical set of parameters, we visualize the 1-parameter family.
  
\section{Derived Schwarz map}\label{sec:DS}
\subsection{Definition of the derived Schwarz map}
\label{sub:defDS}

The equation \ref{eq:Eabc} transforms into 
the SL-form \eqref{eq:SL-form} 
by the projective change of the unknown
\[
    u\longrightarrow \underline{u}:=N\cdot u,\quad\mbox{where}
   \quad N:=\sqrt{x^c(1-x)^{a+b+1-c}}.
\]
The coefficient $q$ is 
\begin{eqnarray*}
  q  &=& -\dfrac14\left\{\dfrac{1-\mu_0^2}{x^2}
         +\dfrac{1-\mu_1^2}{(1-x)^2}
         +\dfrac{1+\mu_\infty^2-\mu_0^2-\mu_1^2}{x(1-x)}\right\}\\
     &=& -\dfrac14\dfrac{(1-\mu_\infty^2)x^2+
           (\mu_\infty^2+\mu_0^2-\mu_1^2-1)x+1
               -\mu_0^2}{x^2(1-x)^2}.
\end{eqnarray*}
Set $v=\underline{u}'$; it satisfies the equation
\begin{equation}\tag*{$dE(a,b,c)$}\label{eq:dEabc}
     v''-\dfrac{q'}qv'-qv=0.
\end{equation}
The Schwarz map of this equation is called the 
{\it derived Schwarz map\/} $DS$ of the equation \ref{eq:Eabc}.

\subsection{Flat fronts and hyperbolic Gauss maps}\label{sub:geom-1}
 Geometrically, the hyperbolic Schwarz map $HS\colon{}X_+\to \HH^3$
 is a flat front in the sense of \cite{KUY},
 and the map $U$ in \eqref{eq:hol-lift} is the {\em holomorphic lift\/}
 of $HS$, which satisfies the differential equation
 \[
       U^{-1}U' = 
                 \begin{pmatrix} 0 & q \\ 1 & 0\end{pmatrix}.
 \]
 That is, under the notations in \cite{KUY}, the
 {\em canonical forms\/} of the flat front $HS$ are
 $\omega=dx$, and $\theta=q\,dx$.
 From now on, we normalize $U$ as $\det U=1$.

 Though $HS$ may have singularities, 
 the {\em unit normal vector field\/}
 $\nu$ is well-defined as
 \begin{equation}\label{eq:normal}
     \nu = U     \begin{pmatrix}
		     1 & \hphantom{-}0 \\
		     0 & -1
		  \end{pmatrix}\transpose{\overline{U}}
 \end{equation}
 The {\em hyperbolic Gauss maps\/} $G$ and $G_*$ are maps defined as
 \begin{align*}
    G(x) & =
      \left(
        \text{
         \begin{tabular}{l}
	  the asymptotic class of the geodesic in $\HH^3$\\
	  starting at $HS(x)$ with initial velocity $\nu(x)$
	 \end{tabular}
        }
      \right)\in \partial H^3=\PP^1\\
    G_*(x) & =
      \left(
        \text{
         \begin{tabular}{l}
	  the asymptotic class of the geodesic in $\HH^3$\\
	  starting at $HS(x)$ with initial velocity $-\nu(x)$
	 \end{tabular}
        }
      \right)\in \partial H^3=\PP^1.
 \end{align*}
 The Schwarz map $S$ and the derived Schwarz map $DS$ are nothing but 
 the hyperbolic Gauss maps $G$ and $G_*$ respectively.
 The ramification points of $DS$ are the zeros of $q$, which
 are the {\em umbilic points\/} of the flat front $HS$.

 The isometric action of $PSL(2,\C)$ to $\HH^3$ as 
 \[
    \HH^3\ni p \longmapsto a p \transpose{{\bar a}}\qquad
    \bigl(\pm a\in PSL(2,\C)\bigr)
 \]
 induces the conformal action on $\partial \HH^3=\PP^1$,
 which coincides with the $PSL(2,\C)$-action as the M\"obius
 transformations.
 Thus, the monodromy representations with respect to 
 $G=S$ and $G_*=HS$ coincide.

\subsection{Description of the image}
\label{sub:image}
Though the maps $S$ and $DS$ are determined by the equation only 
up to linear fractional transformations, in this section, 
we always assume that  
`for any choice of $S$, $DS$ is so chosen that if 
$S=\underline u_1/\underline u_0$ then 
$DS =\underline u'_1/\underline u'_0$'. 

\begin{lemma} Under the convention above
\[
   S(0)=DS(0),\quad S(1)=DS(1),\quad\text{and}\quad
   S(\infty)=DS(\infty).
\]
\end{lemma}
\begin{proof}
 Solutions of the SL-form are singular at $x=0,1,\infty$, 
 so de l'Hopital theorem can be applied. 
\end{proof}

\begin{lemma}
 The equation \ref{eq:dEabc} has the same local behavior with 
 \ref{eq:Eabc} at the three singular points. 
 In addition, it has apparent singularities at the zeros of $q$. 
 Monodromy behaviors of both equations agree. 
 If the zeros are simple, $DS$ ramifies at these points with index $2$; 
 if double, with index $3$.
\end{lemma}
\begin{proof}
In general, the Schwarzian derivative of any Schwarz map $S$
of an equation $w''-Qw=0$ equals $-Q$. From the definition of the 
Schwarzian derivative, it is straightforward to see that $S$ has the
local expression 
\[
S=(x-\xi)^\gamma(\mbox{a non-vanishing holomorphic function at $\xi$})
\]
if and only if $Q$ has the local expression
\[
Q=-\dfrac{1-\gamma^2}{(x-\xi)^2}
+\dfrac{\mbox{a holomorphic function at $\xi$}}{x-\xi}.
\]
 The  SL-form of the equation \ref{eq:dEabc} is given by 
 $\underline v''-\underline q\underline v=0$, where 
 \[
    \underline q 
    =q+\dfrac34\left(\dfrac{q'}q\right)^2-\dfrac12\left(\dfrac{q''}q\right)
    =q+\dfrac14\left(\dfrac{q'}q\right)^2-\dfrac12\left(\dfrac{q'}q\right)'.
 \]
 If $\alpha$ and $\beta$ denote the zeros of $q$, we have
 \[
    \dfrac{q'}{q}=-\frac2x+\frac2{1-x}+\frac1{x-\alpha}+\dfrac1{x-\beta}.
 \]
Since we have
\[
\left(\dfrac{q'}{q}\right)^2=\dfrac4{x^2}+\dfrac{O(1)}x,\quad
\left(\dfrac{q'}{q}\right)'=\dfrac2{x^2}+\dfrac{O(1)}x,
\]
the expression of $\underline q$ above leads to
\[
\lim_{x\to0}x^2q(x)=\lim_{x\to0}x^2\underline q(x)
\quad\left(=-\dfrac{1-\mu_0^2}4\right);
\]
this implies that $S$ and $DS$ have the same local behavior at $0$.
It can be similarly seen that they have the same local behavior 
also at $1$ and $\infty$. When $\alpha\not=\beta$, by a similar 
computation, we have
\[
\lim_{x\to\alpha}(x-\alpha)^2\underline q(x)
 =\lim_{x\to\beta}(x-\beta)^2\underline q(x)
 =\dfrac14+\dfrac12=-\dfrac{1-2^2}4,
\]
and when $\alpha=\beta$,
\[
\lim_{x\to\alpha}(x-\alpha)^2\underline q(x)=-\dfrac{1-3^2}4.
\qquad \mbox{\qedsymbol}
\]
\renewcommand{\qed}{\relax}
\end{proof}

\begin{remark}\label{rem:STW} 
We refer to \cite{STW} for a general treatment of differential
eqautions that admit apparent singularities in addition to
three regular singular points and that the monodoromy groups
are triangle groups.
\end{remark}

The discriminant of the numerator of $q$ is given as
\[
 D:=(\mu_\infty^2+\mu_0^2-\mu_1^2-1)^2-4(1-\mu_0^2)(1-\mu_\infty^2).
\]
This quantity turns out to be symmetric with respect to 
$\{\mu_0^2,\mu_1^2,\mu_\infty^2\}$, and can be expressed as
\begin{multline*}
  D=(s+1)^2-4(t+1),\\
 \text{where}\qquad
  s=\mu_0^2+\mu_1^2+\mu_\infty^2,\quad\text{and}\quad
  t=\mu_0^2\mu_1^2+\mu_1^2\mu_\infty^2+\mu_\infty^2\mu_0^2.
\end{multline*}
Since we assumed $0\le \mu_0^2,\mu_1^2,\mu_\infty^2<1$, 

\begin{figure}
\begin{center}
\includegraphics[width=6cm]{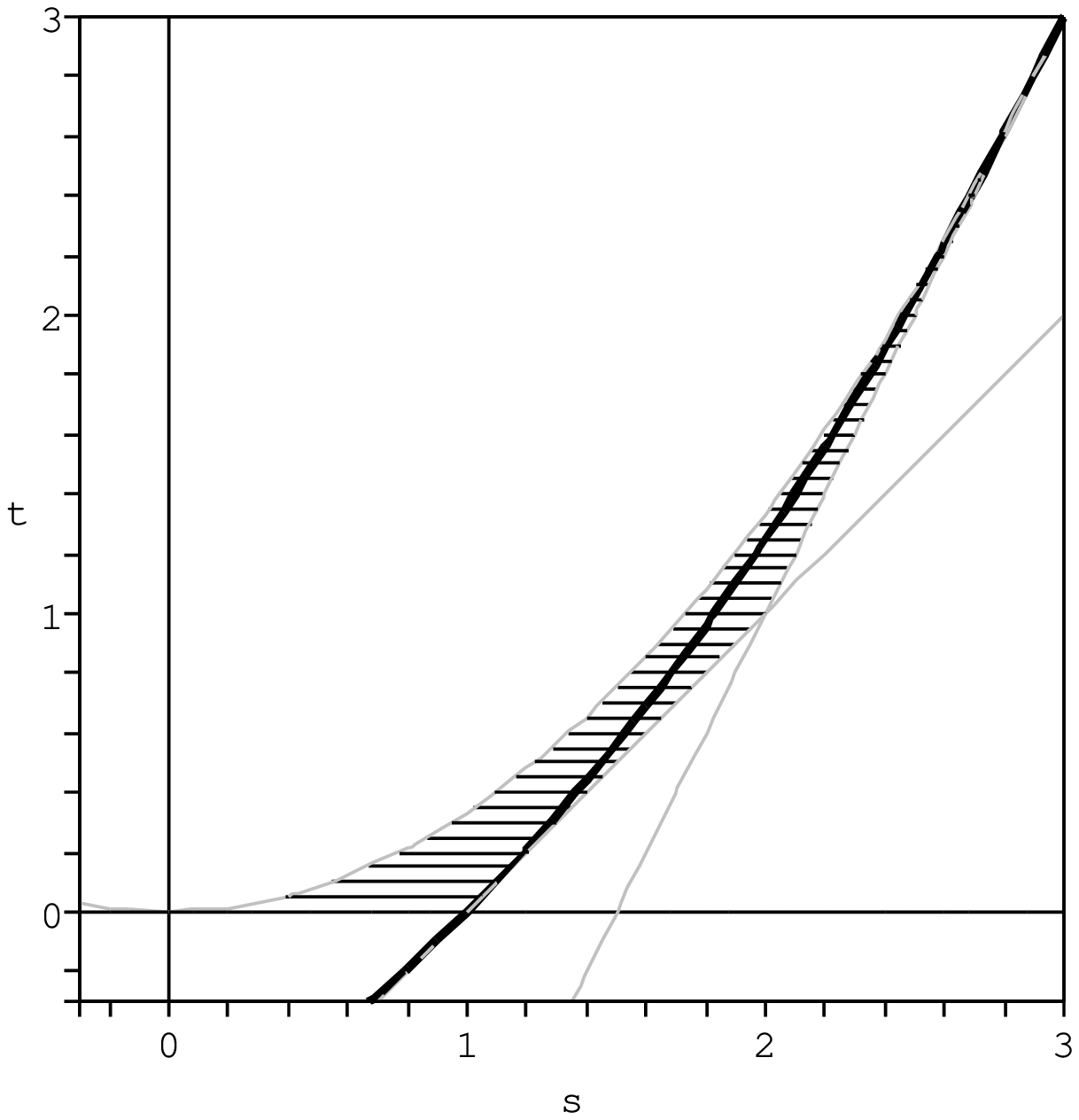}
\hskip1cm
\includegraphics[width=4cm]{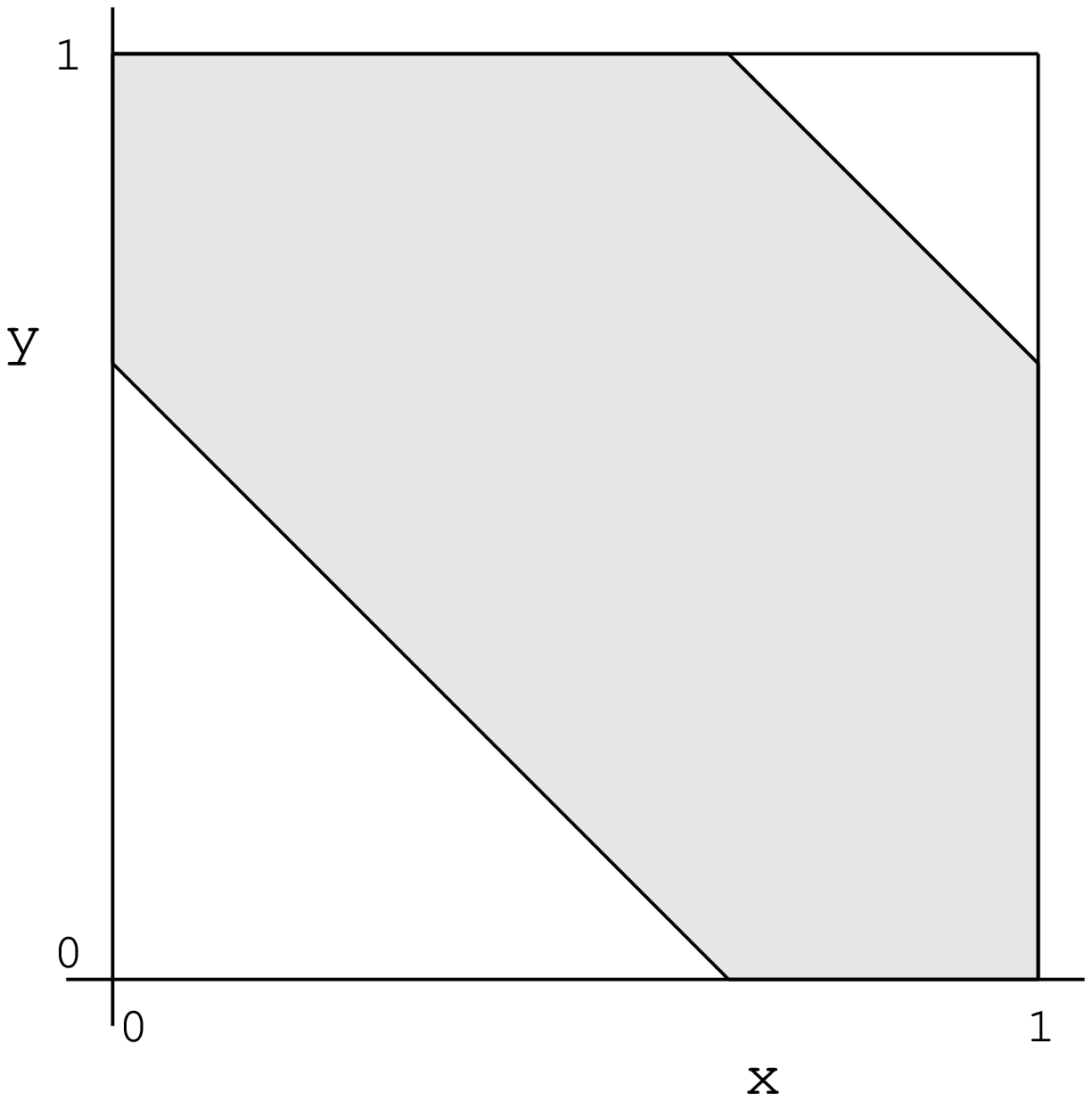}
\end{center}
\caption{Left: Domain of the points $(s,t)$ \qquad Right: Domain of
the points $(x,y)$ for $s=5/3$}
\label{stD}
\end{figure}
\begin{lemma}
 The domain of existence of the point $(s,t)$ is bounded by
 \[
   t=\frac{s^2}3, \quad t=0,\quad t=s-1,\quad\text{and}\quad t=2s-3,
 \]
 see {\rm Figure~\ref{stD}}, 
 and is divided into two parts by the curve $D=0$.
\end{lemma}

\begin{proof}
For notational simplicity, set $s=x+y+z, t=xy+yz+zx$. Substituting
$z=s-(x+y)$ into the second expression, we define the function
\[
t_s(x,y)=xy+(x+y)\{s-(x+y)\},
\]
where $(x,y)$ is in the region:
\[
0\le x,y<1,\quad 0\le s-(x+y)<1.
\]
This region changes its shape at $s=1$ and $2$; so we consider three
cases (1) $0\le s\le1$,  (2) $1\le s\le2$ (refer to the right
picture of Figure~\ref{stD}),
(3) $2\le s<3$. 
In any case,
the function $t_s$ attains its maximum $s^2/3$ at $(s/3,s/3)$, while it
attains its minima $0, s-1, 2s-3$ at different boundary points of the
regions for cases (1), (2), and (3), respectively. 
\end{proof}

We collect properties of the Schwarz and the derived Schwarz maps, 
especially their images.
\begin{proposition}
 For the real parameters $(a,b,c)$ satisfying the condition 
 \eqref{mule1},
 \begin{enumerate}
  \item $S$ gives a conformal isomorphism from $X_+$ onto the triangle
	$T=S(X_+)$,
  \item $DS(x)=S(x)$ for $x=0,1,\infty$,
  \item $DS$, restricted to the interval $(-\infty,0)$ gives a
	diffeomorphsim onto the side $S((-\infty,0))$ of $T$,
  \item $DS$, restricted to the interval $(1,+\infty)$ gives a
	diffeomorphsim onto the side $S((1,+\infty))$ of $T$,
  \item the image $DS((0,1))$ is, as a set, the whole circle $C$
	extending the arc side $S((0,1))$,
  \item the image $DS(X_+)$ is, counting multiplicity, the union of $T$
	and the disc {\rm(}left side of $S((0,1))${\rm )} bounded by $C$,
  \item\label{case:imag} 
	if the set of parameters is in the domain $D<0$, then the
	derived Schwarz map $DS$ has a unique ramification point
	{\rm (}of order $2${\rm )} in the interior of $X_+$, 
	and as $x$ moves from $0$ to $1$, 
	$DS(x)$ moves from $S(0)$ along the circle $C$, 
	in the same sense as $S(x)$ moves from $S(0)$, 
	turns around once the circle and reach at $S(1)$.
  \item\label{case:real}
	if the set of parameters is in the domain $D>0$, then $DS$ has
	two ramification points {\rm (}of order $2${\rm )} on the
	interval $(0,1)$, 
	say $r_1<r_2$, and as $x$ moves from $0$ to $1$, 
	$DS(x)$ moves from $S(0)$ along the circle $C$, 
	in the same sense as $S(x)$ moves from $S(0)$, 
	passes through $S(1)$ and $DS(r_2)$, turns back at $DS(r_1)$ to
	$DS(r_2)$, 
	turns back at $DS(r_2)$ to $S(0)$ and eventually ends the
	journey at $S(1)$,
  \item if the set of parameters satisfy $D=0$, 
	then $r_1=r_2$ is the unique ramification point 
	{\rm (}of order $3${\rm )}
	in $X_+$; 
	this is the limit case of \ref{case:imag} and \ref{case:real}.
\end{enumerate}
\end{proposition}
\begin{proof}
The assertion (1) is well-known and (2) is Lemma 2.1. 
 We prove $(3)$--$(9)$.
 We first assume $D<0$, and study the images of one of the three 
 intervals, say $(0,1)$ under $S$ and $DS$. 
 Since the local behavior (and the values) of 
 $S$ and $DS$ coincide at each $x=0,1,\infty$, 
 and since holomorphic maps preserves orientation, 
 as $x$ moves from $0$ to $1$, $DS(x)$ moves from $S(0)$ to $S(1)$ 
 along the circle $C$ in the same sense as $S(x)$ moves from $S(0)$ to $S(1)$. 
 The point $DS(x)$, after leaving $S(0)$ either stops at $S(1)$, 
 or passes through $S(1)$ and turns around the whole circle, 
 passes through $S(0)$ and then stops at $S(1)$, 
 or turns once more, \dots.
 The same assertion holds for the other intervals $(-\infty,0)$ 
 and $(1,+\infty)$. 

 \noindent
 On the other hand, since there is only one ramification point of degree 2 
 in $X_+$, the rotation number is 2. 
 Thus two out of the three intervals are mapped by $DS$ to the arc, 
 the image under $S$, and only one interval is mapped to the arc 
 plus the whole circle. 
 In the special dihedral ($n=3$) case 
 where $(\mu_0,\mu_1,\mu_\infty)=(1/2,1/2,1/3)$,
 it is known that the interval $(0,1)$ 
 plays this part;  refer to \cite{SYY}.
 Since the domain $D<0$ is connected, 
 we conclude that it is always the case.

\noindent
 We next consider what happens when the zeros $\alpha\in \C-\R$ and 
 $\beta=\bar\alpha$ come together ($D=0$), and then turn 
 to two real roots ($D<0$). 
 During this process, global behavior of $DS$ does not change. 
 Local behavior around the zeros can be best understood 
 via the following model, by which the proposition is readily proved. 
\end{proof}

\subsection{A model of confluence of the two ramification 
points of the derived Schwarz map} \label{sec:model}
\begin{figure}
\footnotesize
\begin{center}
\begin{tabular}{ccc} 
 \includegraphics[width=5cm]{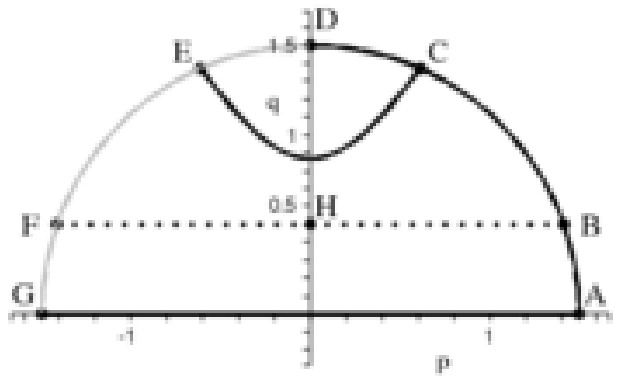} & 
  \raisebox{25mm}{$\rightarrow$} &
 \includegraphics[width=5cm]{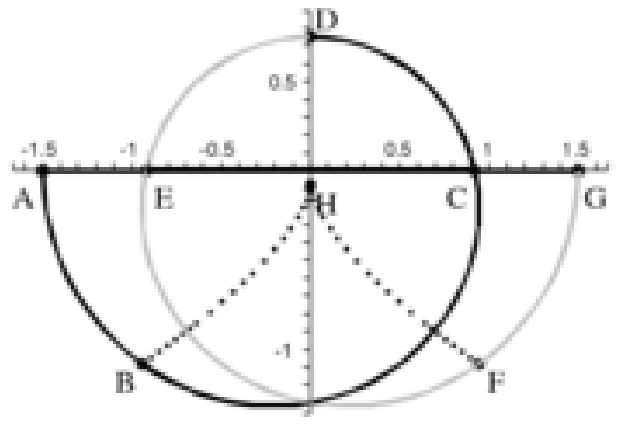} 
\end{tabular}
\end{center}
\caption{The image of hemi-disc by $x\mapsto z=-x^3/3 - x/4$}
\label{DS.conf1}
\end{figure}

\begin{figure}
\footnotesize
\begin{center}
\begin{tabular}{ccc} 
 \includegraphics[width=5cm]{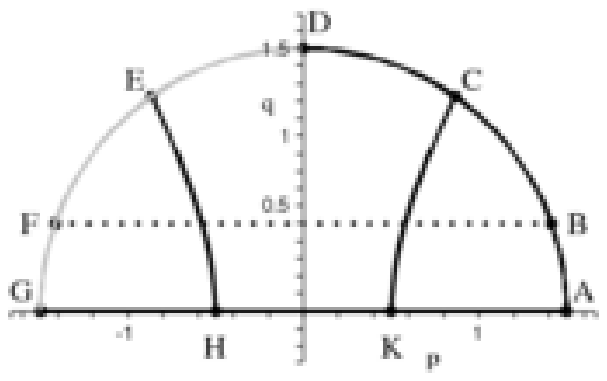} & 
  \raisebox{25mm}{$\rightarrow$} &
 \includegraphics[width=5cm]{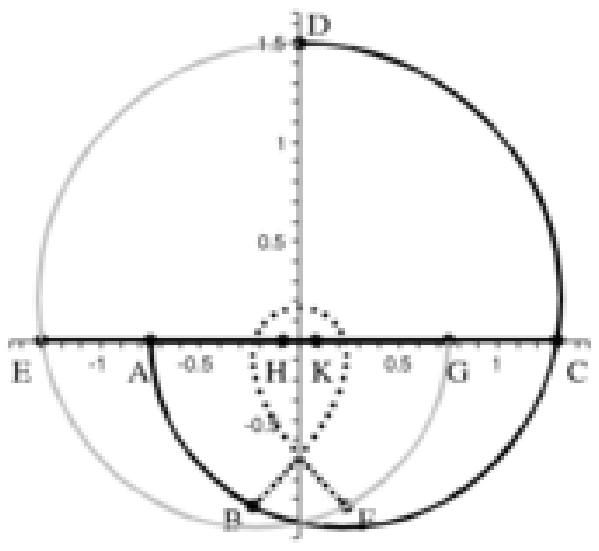} 
\end{tabular}
\end{center}
\caption{The image of hemi-disc by $x\mapsto z=-x^3/3 + x/4$}
\label{DS.conf2}
\end{figure}

We see what will happen for the derived Schwarz map when the two
complex conjugate ramification points come together and separate 
into two real points. Such a map (a family of maps) can be locally 
expressed by
\[
\phi_t: \C\ni x\mapsto z=-\left(\dfrac{x^3}3+tx\right)\in \C,
\]
where the real parameter $t$ varies from positive to negative.
In fact since $z'=-(x^2+t)$, 
the ramification points $\pm i\sqrt t$ come 
together and then separate to $\pm\sqrt{-t}$.  
Figures \ref{DS.conf1} and \ref{DS.conf2} show the images 
of the upper unit hemi-disc when $t=1/4$ and $t=-1/4$.
Points marked alphabets are mapped to the points with the same marks.

In Figure \ref{DS.conf1}, 
the curve $\stackrel{\frown}{EC}$ in the left figure 
is the preimage of the segment $\overline{EC}$ on the real axis
in the right figure.
It divides 
the hemi-disc into two regions; the upper one covers once 
the upper part of the $z$-plane, 
while the lower one covers twice the lower part.
$H$ denotes the ramification point.

\par\noindent
In Figure \ref{DS.conf2}, 
the curves $\stackrel{\frown}{EH}$ and $\stackrel{\frown}{CK}$ 
in the left figure are the preimages of
the line segments $\overline{EH}$ and $\overline{CK}$
on the real axis, respectively. 
They divide the hemi-disc into three regions; 
the middle one covers the upper part of the $z$-plane once, 
and each of the left and the right ones covers the lower part once.
The point $H$ and $K$ are the ramification points.

As parameters vary and pass through the curve $\{D=0\}$, 
near the two confluenting ramification points, locally, 
the derived Schwarz map behaves similar to $\phi_t$.

\subsection{Illustration of the image of the derived Schwarz map}
We illustrate the behavior of the derived Schwarz map when the set of
parameters traverses $D=0$. 
We study the (derived) Schwarz map with 
parameters $(a,b,c)=(1/2,1/2,c)$; we let $c$ move from 1 to 0.  
In the $st$-plane in Figure~\ref{stD}, 
this move corresponds to the path along $t=s^2/4$ 
starting at $(0,0)$; it crosses the curve $D=0$ at 
$(3/2,9/16)$ ($c=1-\sqrt3/2$), and goes out of the domain at $(2,1)$.

\subsubsection{A normalization of the {\rm (}derived{\rm )} Schwarz map}
\label{subsub:normalization}
In this subsection, for generic parameters $(a,b,c)$, 
we {\it normalize} the Schwarz map as
\[
       S(0)=0,\quad S(1)=1,\quad \text{and}S(\infty)=\infty,
\]
The Schwarz map, in general,
is multi-valued and determined by the equation 
\ref{eq:Eabc} only up to conjugacy. 
Here we first choose two solutions
\[
    u_0:=F(a,b,c,x) \quad\mbox{and}\quad 
    u_1:=x^{1-c}F(a-c+1,b-c+1,2-c,x)
\]
in the upper half $x$-plane 
$X_+:=\{x\in\C\mid \Im(x)\ge0, x\not=0,1\}$,
where $F(a,b,c,x)$ denotes the hypergeometric function
with parameters $(a,b,c)$.
They are real valued on the interval $(0,1)$, since the parameters
are real.
We define temporarily a Schwarz map: $S_T(x):=u_1/u_0(x)$. 
The image of the interval $(0,1)$ is a real interval,
$S_T(0)=0$, and the image of the interval $(-\infty,0)$ 
is a segment.
\begin{lemma} 
 If $a+b-c>0$, then
 \[
     S_T(1)=\dfrac{\Gamma(2-c)\Gamma(a)\Gamma(b)}
        {\Gamma(c)\Gamma(a-c+1)\Gamma(b-c+1)}.
 \]
 If $0<c<2$ and $a\le b$, then
 \[
      S_T(\infty)=\dfrac{\Gamma(2-c)\Gamma(b)\Gamma(c-a)}
          {\Gamma(c)\Gamma(1-a)\Gamma(1-c+b)}e^{\pi i (1-c)}.
 \]
\end{lemma}
\begin{remark}
 Though the left hand-sides of the formulae are, 
 by definition, symmetric in $a$ and $b$, 
 the right hand-side of the second formula not.
\end{remark}
\begin{proof}
(1) Consider solutions around $x=1$:
 \begin{align*}
     v_0&:=F(a,b,a+b-c+1;1-x)\quad\text{and}\\\
     v_1&:=(1-x)^{c-a-b}F(c-a,c-b,c+1-a-b;1-x);
 \end{align*}
 set $s_0=u_1/u_0, s_1=v_1/v_0$.
 They are related as
 \[
    s_0 = \dfrac{A s_1+C}{B s_1+D}\quad\mbox{along} (0,1),
 \]
 where
 \begin{alignat*}{2}
  D&:=\dfrac{\Gamma(c)\Gamma(c-a-b)}{\Gamma(c-a)\Gamma(c-b)},\qquad
  &C&:=\dfrac{\Gamma(2-c)\Gamma(c-a-b)}{\Gamma(1-a)\Gamma(1-b)}\\
  B&:=\dfrac{\Gamma(c)\Gamma(a+b-c)}{\Gamma(a)\Gamma(b)},\qquad
  &A&:=\dfrac{\Gamma(2-c)\Gamma(a+b-c)}{\Gamma(a-c+1)\Gamma(b-c+1)}.
 \end{alignat*}
 If $a+b-c>0$, then $s_1(1)=\infty$. 
 Thus we have $s_0(1)=A/B$. 

 (2) We make use of Kummer's relations:
 \begin{align*}
   u_0&=(1-x)^{-a}F(c-b,a,c;\dfrac{x}{x-1}),\\
  u_1&=x^{1-c}(1-x)^{c-a-1}F(1-b,a+1-c,2-c;\dfrac{x}{x-1}),
 \end{align*}
 and let $x\to-\infty$. 
 Then we have
 \[
    \dfrac{x^{1-c}(1-x)^{c-a-1}}{(1-x)^{-a}}\longrightarrow 
    e^{\pi i(1-c)} 
    \quad\mbox{and} \quad \dfrac{x}{x-1}\nearrow 1.
 \]
 The Gauss formula
 \[
     F(a,b,c;1)=\dfrac{\Gamma(c)\Gamma(c-a-b)}{\Gamma(c-a)\Gamma(c-b)},
    \quad\mbox{if}\quad c>0,\ c-a-b>0
 \]
 tells that, under the conditions $0<c<2$ and 
 \[
    c-(c-b)-a=(2-c)-(1-b)-(a+1-c)=b-a>0,
 \]
 we have
 \begin{align*}
    \dfrac{F(1-b,a+1-c,2-c;1)}{F(c-b,a,c;1)}
    &=\dfrac{\Gamma(2-c)\Gamma(b-a)}{\Gamma(b-c+1)\Gamma(1-a)}\cdot
      \dfrac{\Gamma(b)\Gamma(c-a)}{\Gamma(c)\Gamma(b-a)}\\[3mm]
    &=
     \dfrac{\Gamma(2-c)\Gamma(b)\Gamma(c-a)}
      {\Gamma(c)\Gamma(b-c+1)\Gamma(1-a)}.
      \qquad\mbox{\qedsymbol}
 \end{align*}
\renewcommand{\qed}{\relax}
\end{proof}
We then define a new Schwarz map $S$ which sends $0$, $1$ and 
$\infty$ to $0$, $1$ and $\infty$, respectively, by
\[
    S(x)=\dfrac{S_T(x)}{S_T(x)-S_T(\infty)}
      \dfrac{S_T(1)-S_T(\infty)}{S_T(1)}.
\]
Accordingly, the corresponding new derived Schwarz map $DS$ is defined:
put
\[
   U_0:=N\cdot u_0,\quad U_1:=N\cdot u_1,
       \quad DU_0:=(U_0)',\quad DU_1:=(U_1)',
\]
and $DS_T:=DU_1/DU_0,$  and we define 
\[
    DS(x)=\dfrac{DS_T(x)}{DS_T(x)-S_T(\infty)}
       \dfrac{S_T(1)-S_T(\infty)}{S_T(1)}. 
\]

\subsubsection{Illustration of confluence}
\label{subsub:confluence}
\begin{figure}
\begin{center}
\includegraphics[width=8cm]{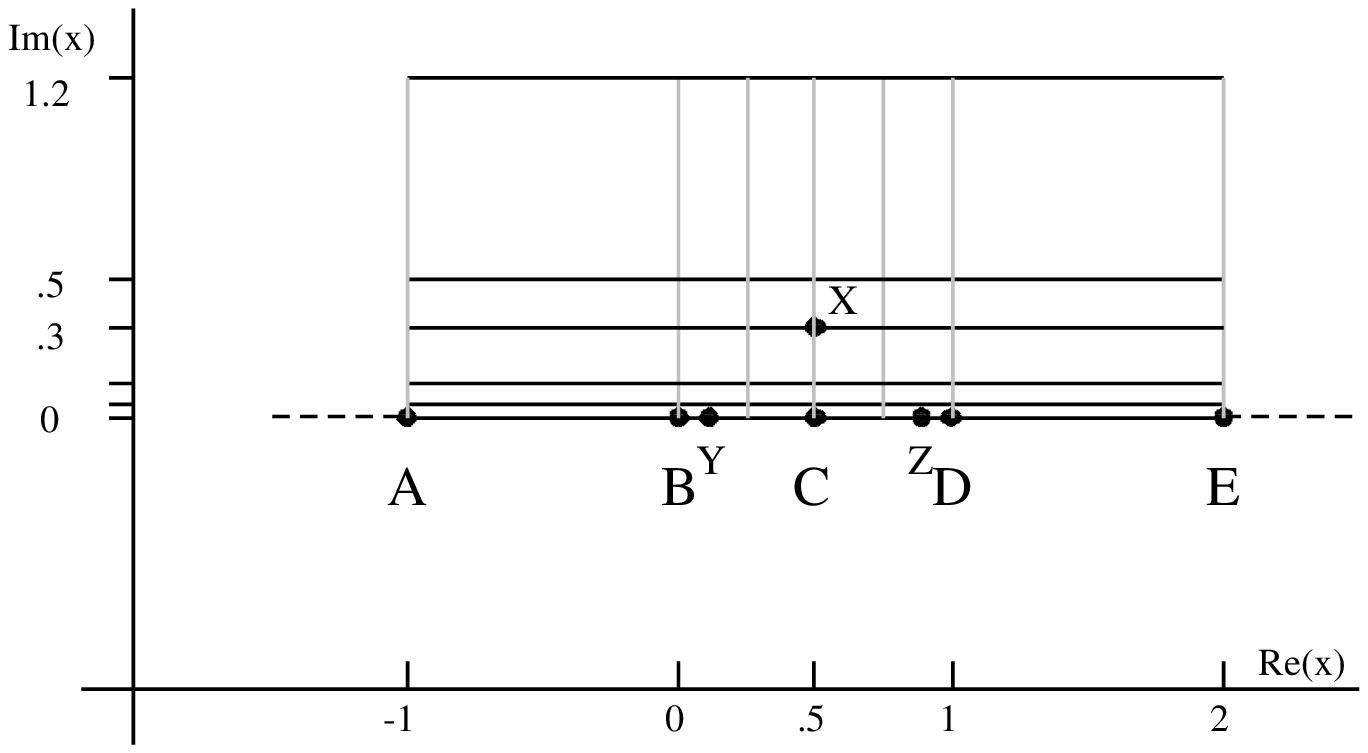}
\end{center}
\caption{}
\label{upper}
\end{figure}

We illustrate the images of the upper half plane under $S$ and $DS$ with
parameters $(a,b,c)=(1/2,1/2,c)$, $c=1,\dots,0.05$.
Note that in this case,
\[
   \mu_0^2=\mu_1^2=(1-c)^2\quad\text{and}\quad \mu_\infty^2=0.
\]
Since $q=0$ is written as $x^2-x+(2c-c^2)=0$, 
the discriminant is $D=4(1-c)^2-3$; we thus have
\[
   D|_{1>c\ge1-\sqrt3/2}\le0,\qquad D|_{1-\sqrt3/2\ge c>0}\ge0.
\]
Namely, when $c\ge 1-\sqrt{3}/2$, the ramification point is
located at 
\[
   \frac{1}{2}\left(1+i\sqrt{-4c^2+8c-1}\right),
\]
and when $c< 1-\sqrt{3}/2$, we have two ramification points
\[
    \frac{1}{2}\left(1\pm\sqrt{4c^2-8c+1}\right)
\]
on the real axis.
Let us take a domain in the upper half plane as in Figure \ref{upper},
where $A=(-1,0)$, $B=(0,0)$, $C=(1/2,0)$, $D=(1,0)$ and $E=(2,0)$, and
the height of the quadrangle is $10/8$. 

\begin{figure}
\footnotesize
\begin{center}
\begin{tabular}{c|cc} 
 value of $c$   & image under $S$ & image under $DS$\\ \hline
 \raisebox{20mm}{0.90}
                & \includegraphics[width=5cm]{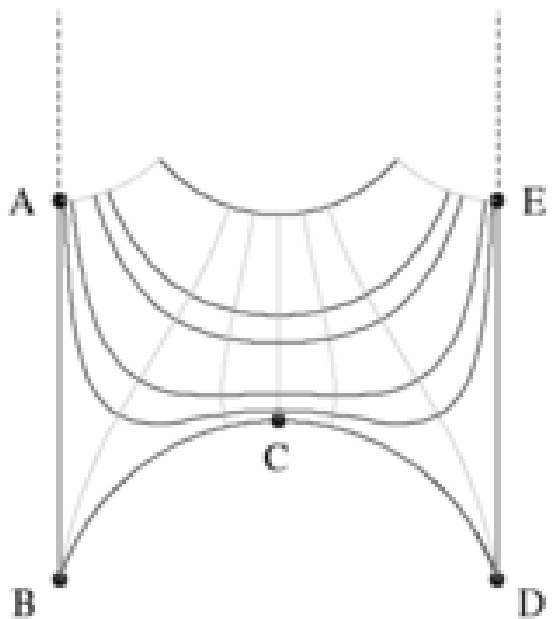} 
                & \includegraphics[width=5cm]{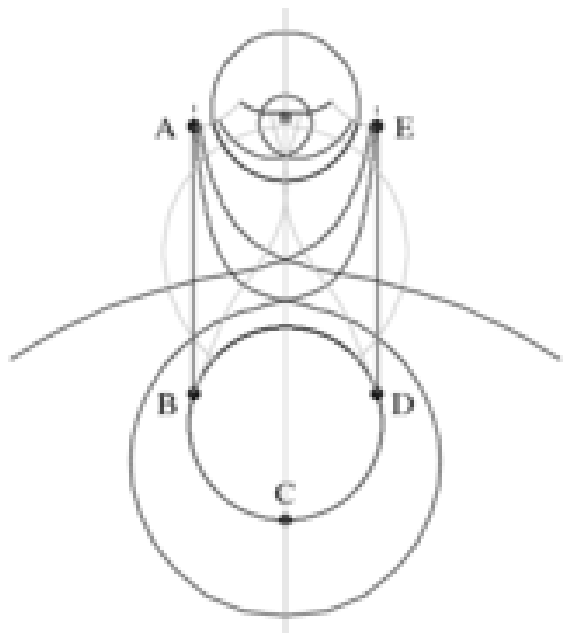} \\
\hline
 \raisebox{20mm}{0.50}
                & \includegraphics[width=5cm]{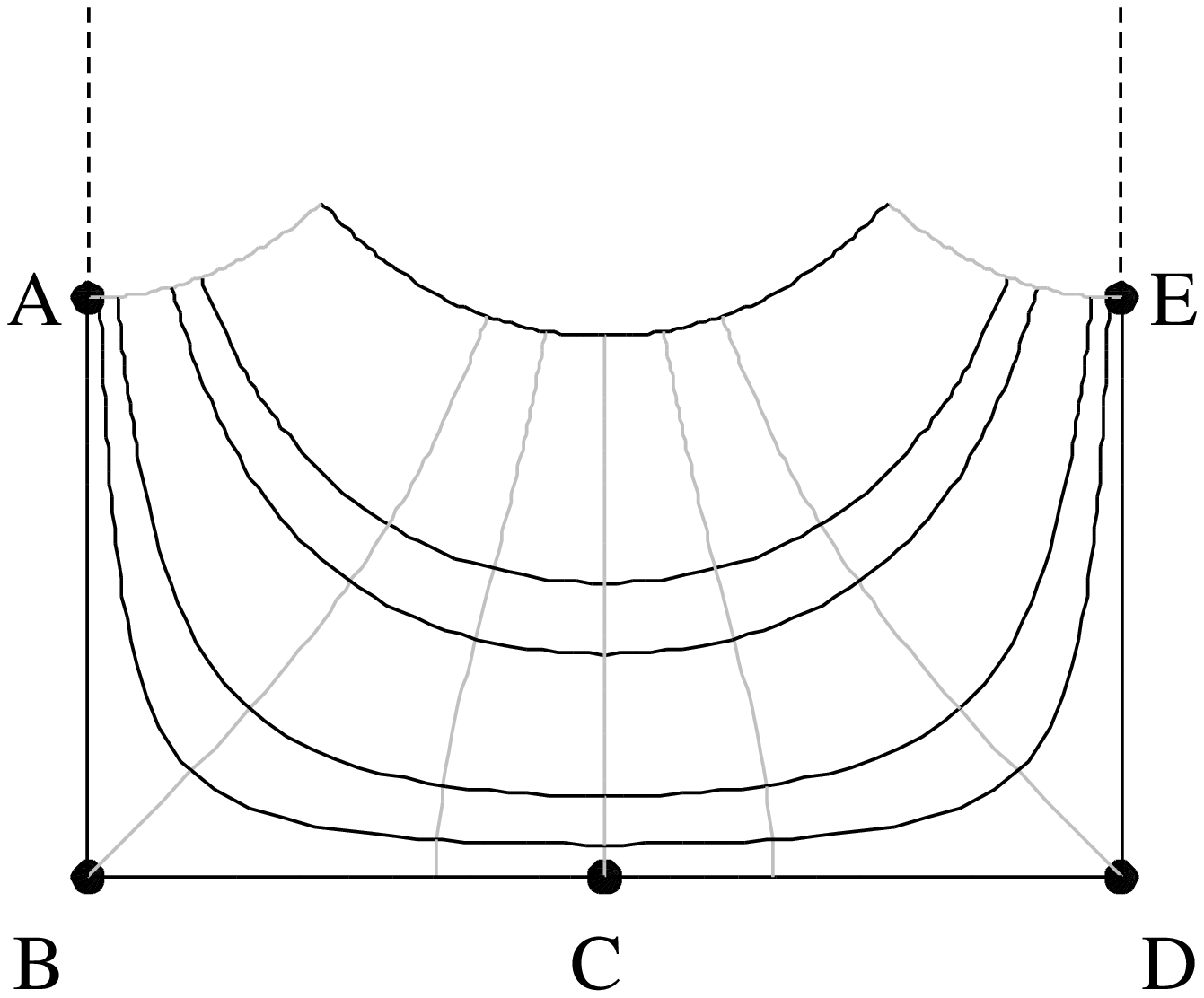}
                & \includegraphics[width=5cm]{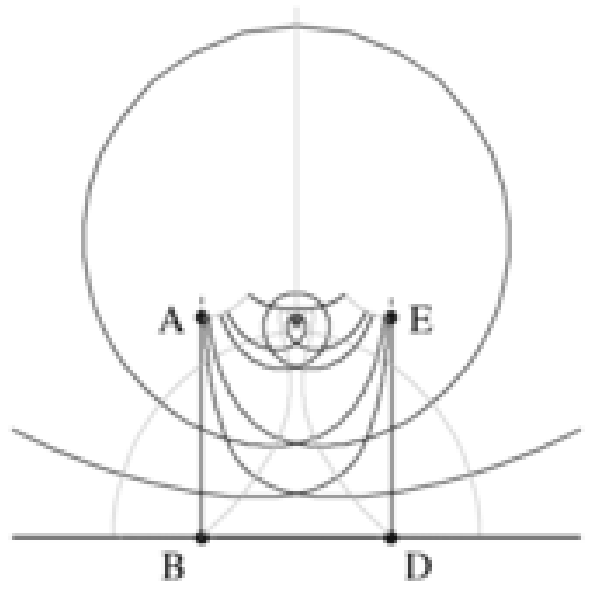} \\
\hline
 \raisebox{20mm}{0.20}
                & \includegraphics[width=5cm]{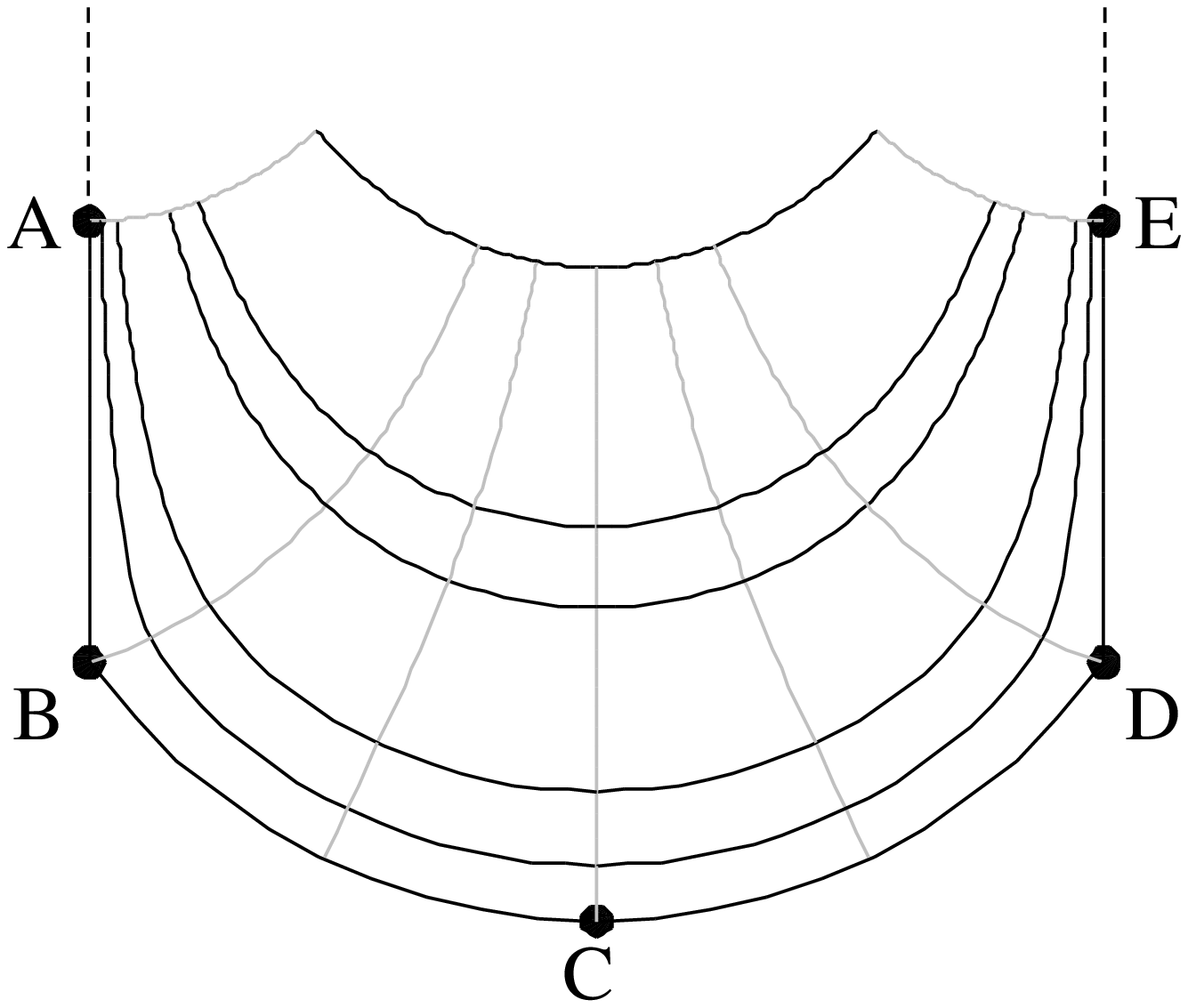} 
                & \includegraphics[width=5cm]{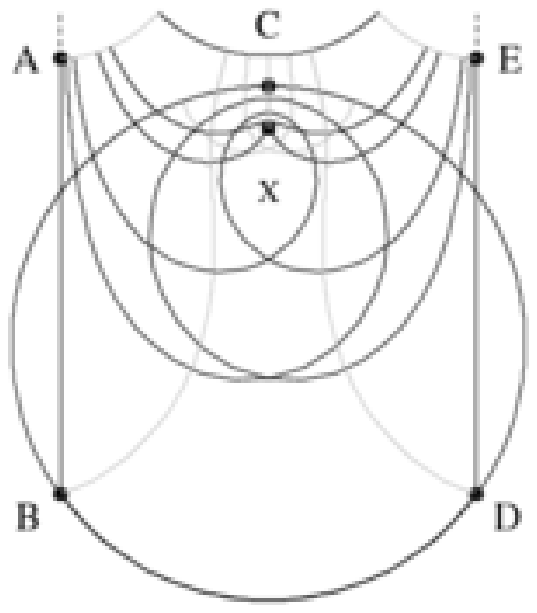} \\
\hline
\end{tabular}
\end{center}
\caption{Images of the domain in Figure \ref{upper} under $S$ and $DS$
when $c>1-\sqrt{3}/2$}
\label{illst1}
\end{figure}

\begin{figure}
\footnotesize
\begin{center}
\begin{tabular}{c|cc} 
 value of $c$   & image under $S$ & image under $DS$\\ \hline
 \raisebox{20mm}{0.05}
                & \includegraphics[width=5cm]{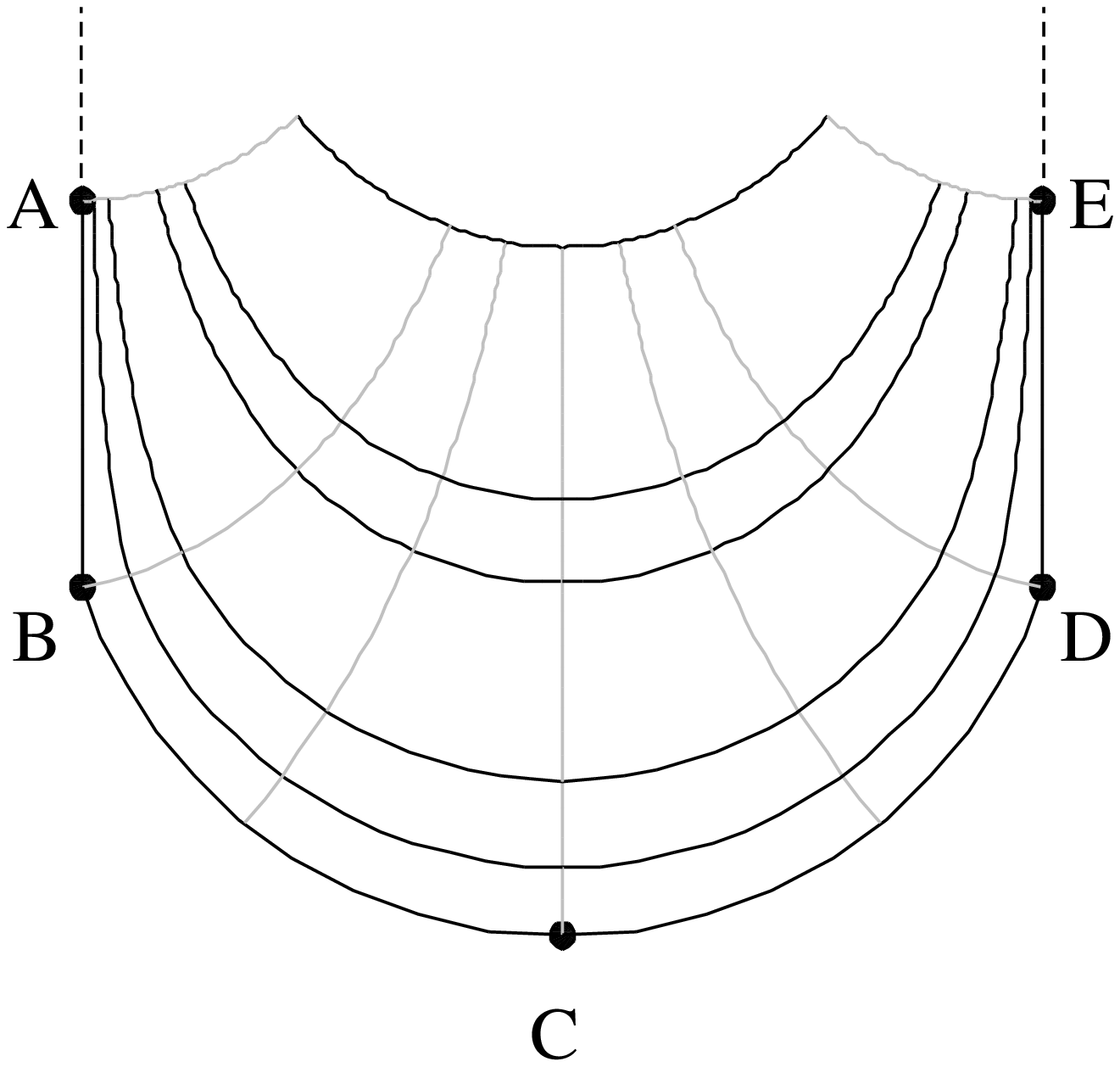} 
                & \includegraphics[width=5cm]{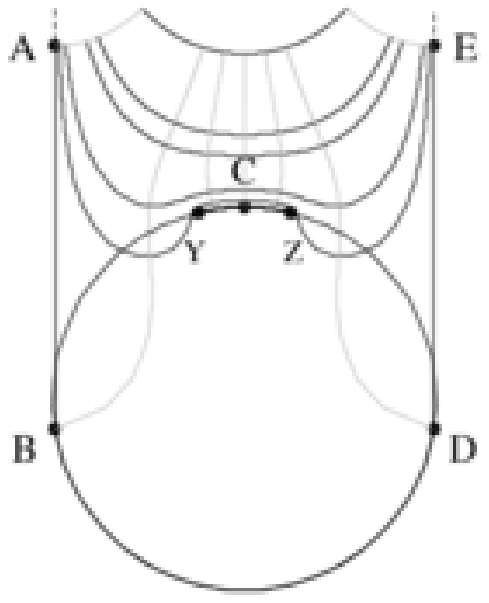} \\
\hline
\end{tabular}
\end{center}
\caption{Images of the domain in Figure \ref{upper} under $S$ and $DS$
when $c=0.05 (<1-\sqrt{3}/2)$}
\label{illst2}
\end{figure}

The point $X=(0.5,\sqrt{0.11})\sim(0.5, 0.3317)$ 
denotes the ramification point when $c=0.2$ 
and the points $Y=(0.5-\sqrt{0.1525},0)\sim(0.1095,0)$ and 
$Z=(0.5+\sqrt{0.1525},0)\sim(0.8905,0)$ denote
the ramification points when $c=0.05$. 
The images of these points under the derived Schwarz maps 
are bullets in Figures \ref{illst1} and \ref{illst2}
with the same name. 
The bullets in Figure \ref{illst1} ($c=0.90,\ 0.50$) without names are the
(images of) ramification points.

%

Figures \ref{illst1} and \ref{illst2} show
how the images of the domain under $S$ and $DS$ 
depend on the parameter $c$.
In Figure \ref{illst1} the image of the segments $\overline{AB}$,
$\overline{BC}$, $\overline{CD}$ and $\overline{DE}$
under $DS$ is the segment $\overline{AB}$, the arc
${BDC}$, the arc ${CBD}$ 
and the segment $\overline{DE}$, respectively.
However, in Figure \ref{illst2} when $c=0.05$, the images of $\overline{BC}$
and $\overline{CD}$ are the curves
\[
BDZCYC\quad\mbox{and}\quad CZCYBD.
\]
In these figures, around $C$ (and $X,Y,Z$), we can observe the happenings 
described in Section \ref{sec:model} using the model map $\phi_t$.

In the case $c=1$, the inverse of the Schwarz map $S$ is the elliptic
modular function $\lambda$, of which behavior is well-known.
But the Schwarz map itself is a bit difficult to treat. 
The solutions $u_0$ and $u_1$ coincide, and the solutions $v_0$ and
$v_1$ coincide.
So the connection matrix relating these solutions loses its sense.
Instead, we draw the associated surfaces in the $3$-ball as
follows.

\begin{figure}
\begin{center}
\includegraphics[width=5.5cm]{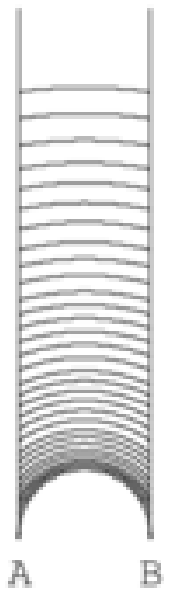}
\end{center}
\caption{The image $S(X_+)$}
\label{lam1}
\end{figure}
\begin{figure}
\begin{center}
\includegraphics[width=5.5cm]{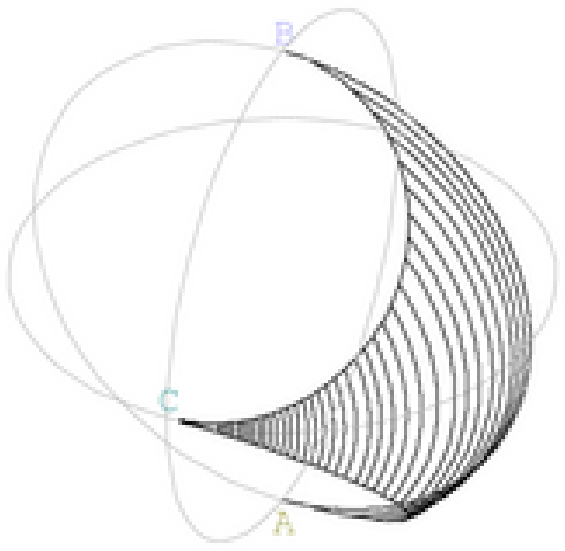}
\hskip0.5cm
\includegraphics[width=5.5cm]{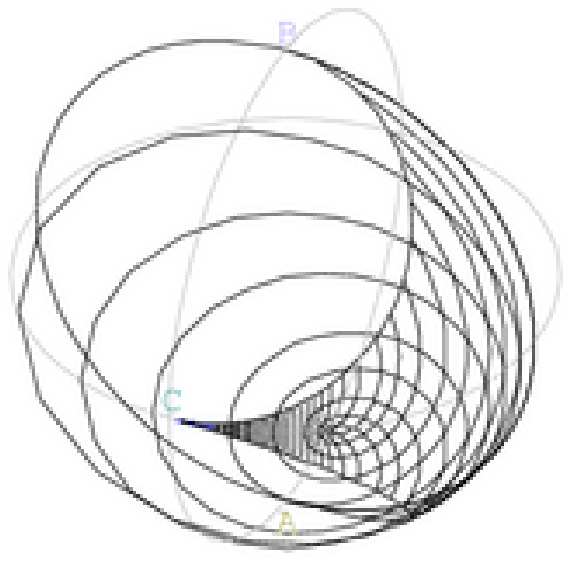}
\end{center}
\caption{The images $\chi(S(X_+))$ and $\chi(DS(X_+))$}
\label{lam2}
\end{figure}

In order to get a total view of the images, we identify $\PP^1$ with the
boundary of $3$-ball by the mapping
\[
    \chi: \PP^1\ni z \longmapsto 
    \left({z+\overline{z}\over 1+z\overline{z}},
      {-i(z-\overline{z})\over 1+z\overline{z}},
       {z\overline{z}-1\over 1+z\overline{z}}\right)
      \in \partial\BB^3, 
\]
which is compatible with identification of $\HH^3$ with 
the ball $\BB^3$.
In the case where $(a,b,c)=(1/2,1/2,1)$, the image of
the upper half plane under $S$ is the fundamental domain relative
to the triangular group usually denoted by $(\infty,\infty,\infty)$
that has the picture in Figure \ref{lam1}.
Its image in $\partial{\bf B}^3$ looks as in the left figure 
of Figure \ref{lam2}.
The image under $DS$ realized in $\partial{\bf B}^3$
is as in the right of Figure \ref{lam2}.
Here, the three gray-colored great circles on 
$\partial{\bf B}^3$ are added to evoke the reader a stereographic image.

\section{%
 Parallel family of flat fronts\\
 connecting Schwarz and derived Schwarz maps
}
In this section, we study relations among
the Schwarz map, the derived Schwarz map and 
the hyperbolic Schwarz map.

\subsection{A relation between $S$ and $DS$}
Since $S$ restricted to $X_+$ is a biholomorphic isomorphism 
between $X_+$ and a Schwarz triangle $T=S(X_+)$, 
we can consider the composite map 
$f:=DS\circ S^{-1}|_{X_+}$. 
This map has a simple expression
\begin{equation}\label{eq:s-and-ds}
   f: T\ni z\longmapsto z+2\dfrac{\dot{x}}{\ddot{x}},
\end{equation}
where $\dot{~} = d/dz$. 
To show this expression, we let
$u$ and $v$ be solutions of an $SL$-equation \eqref{eq:SL-form}
such that $uv'-uv'=1$. 
Since $z'(:=dz/dx)=-1/v^2$ and $\ddot{x}=d^2x/dz^2$, we have 
$v=i/{\sqrt z'}=i\sqrt{\dot{x}}$, $u=vz$, 
and 
\[
  u'=i\dfrac1{\sqrt{\dot{x}}}+z\dfrac{i}2(\dot{x})^{-3/2}\ddot{x},
  \qquad
  v'=\dfrac{dv}{dx}=
  \dfrac{dv}{dz}\dfrac{dz}{dx}
  =\dfrac{i}2(\dot{x})^{-3/2}\ddot{x}.
\]
Hence, we get \eqref{eq:s-and-ds}.
\begin{proposition}
 If $x(z)$ is invariant under $g\in PSL(2,\C)$, that is, 
 $x(gz)=x(z)$, then $f(z)$ is covariant under $g${\rm :}
 $f\circ g(z)=g\circ f(z)$.
\end{proposition}
\begin{proof}
 Set 
 \[
    g(z)=\dfrac{az+b}{cz+d},\qquad ad-bc=1.
 \]
 Then we have
 \[
           \dot{x}(gz)=(cz+d)^2\dot{x}(z),\quad
    \ddot{x}(gz)=(cz+d)^4\ddot{x}(z)+2c(cz+d)^3\dot{x}(z).
 \]
 Compute $f(gz)$ using these formula, we have the conclusion.
\end{proof}

\begin{remark}
 The covariance follows from the way of constructing the
 hyperbolic Gauss map, once we know that $S$ and $DS$
 are the hyperbolic Gauss maps of the hyperbolic Schwarz map $HS$
 (see Section~\ref{sub:geom-1}).
\end{remark}

\subsection{Parallel family of flat fronts}
\label{sub:flat-fronts}
As seen in Section~\ref{sub:geom-1}, the hyperbolic Schwarz map
is considered as a flat front
\[
   \varphi:=HS=U\transpose{\overline{U}}\colon{}X\longrightarrow \HH^3
\]
where $U\colon{}X\to SL(2,\C)$ is the holomorphic lift 
as in \eqref{eq:hol-lift}.

The {\em parallel front\/} $\varphi_t$ of distance $t\in\R$  of
the flat front $\varphi$ is defined as
\begin{equation}\label{eq:parallel}
   \varphi_t(x) = \exp_{\varphi(x)}t\nu(x)
                = (\cosh t) \varphi(x) + (\sinh t)\nu(x),
\end{equation}
where $\nu$ is the unit normal vector field
as in \eqref{eq:normal} and $\exp$ denotes the exponential map of
$\HH^3$.
In the second expression of \eqref{eq:parallel}, 
we identify $\HH^3$ with the upper half component of the 
two-sheet hyperboloid in the Minkowski $4$-space.
Then the parallel front is expressed as
\[
    \varphi_t = U_t\transpose{\overline{U_t}},\qquad
    U_t = U\begin{pmatrix}
	      e^{t/2}  & 0 \\ 0 & e^{-t/2}
	   \end{pmatrix},
\]
see \cite{KUY}, in which the holomorphic lift is denoted by 
``$E$''.
The hyperbolic Gauss maps 
of $\varphi_t$ coincide with those 
of $\varphi$ (which are $G=S$ and $G_*=DS$ in our case),
and the normal geodesic starting at $\varphi(x)$ intersect with 
the ideal boundary at $G(x)=S(x)$ and $G_*(x)=DS(x)$.

The induced metric (the first fundamental form) of $\varphi_t$
is expressed as
\[
    ds^2_t = q_t\,dx^2 + \bar q_t\,d\bar x^2 +
              (1+q_t\bar q_t)\,dx\,d\bar x,\qquad
	      \text{where}\qquad
	      q_t = e^{t}q.
\]
Set $r(x):=-\log|q|$.
Then $\varphi_{r(x)}$ has a singularity at $x$.
Hence the locus of the singular points of $\varphi_t$ is expressed as
\[
   \psi(x) = \varphi_{r(x)}(x),
\]
which is the {\em caustic\/} of the front $\varphi$.
It is known that, locally, 
the caustic of the flat front is a flat front,
see \cite{R}.
More detailed discussions are found in \cite{KRUY};
refer also to \cite{GMM, KUY, KRSUY} for the materials above.

\begin{figure}[t]
\begin{center}
\includegraphics[width=3.6cm]{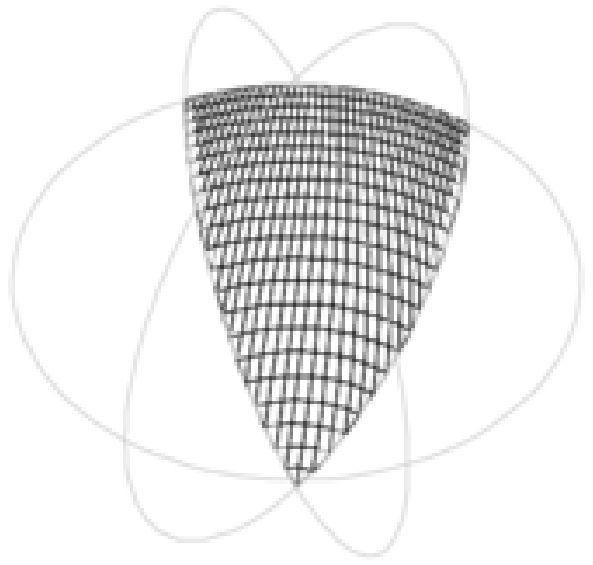}
\hskip.2cm
\includegraphics[width=3.6cm]{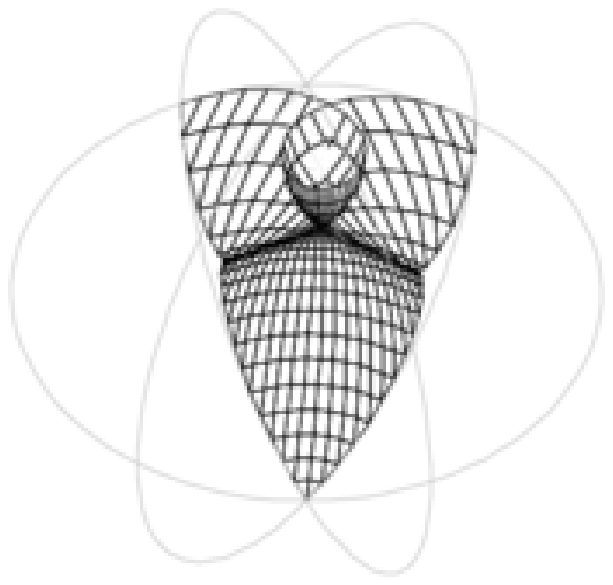}
\hskip.2cm
\includegraphics[width=3.6cm]{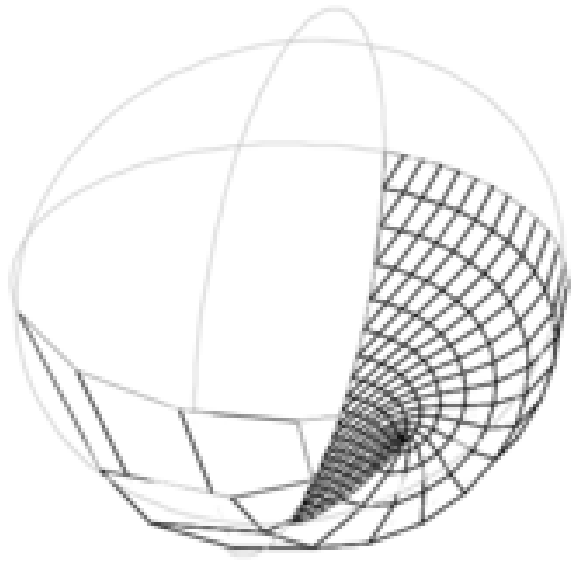}
\end{center}
\caption{The images $\chi(S(X_+))$, $HS(X_+)$ and $\chi(DS(X_+))$}
\label{dihed1}
\end{figure}

\begin{figure}
\footnotesize
\begin{center}
\begin{tabular}{ccc} 
 \includegraphics[width=3.6cm]{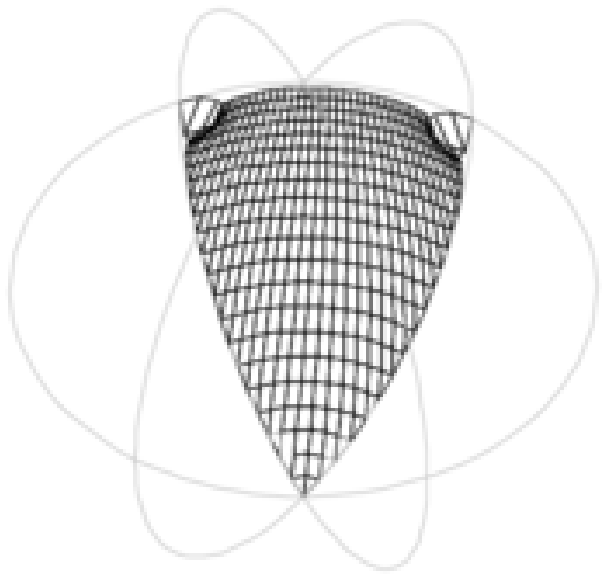} & 
 \includegraphics[width=3.6cm]{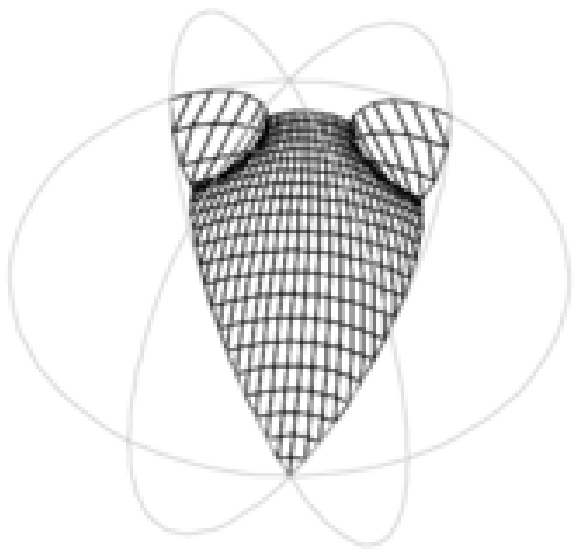} &
 \includegraphics[width=3.6cm]{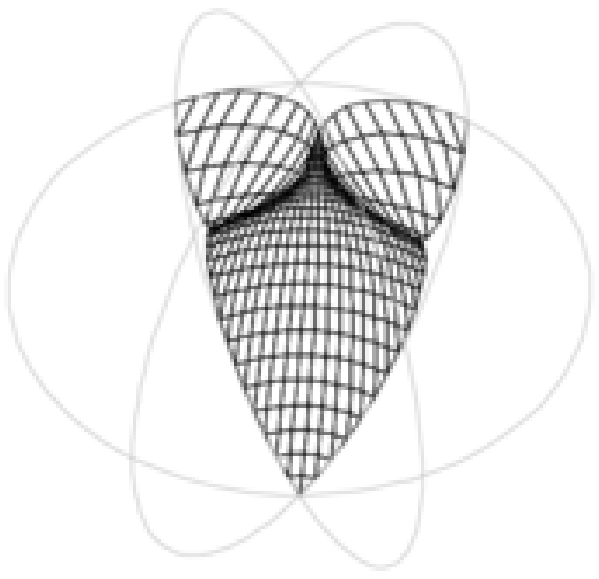} \\
 \includegraphics[width=3.6cm]{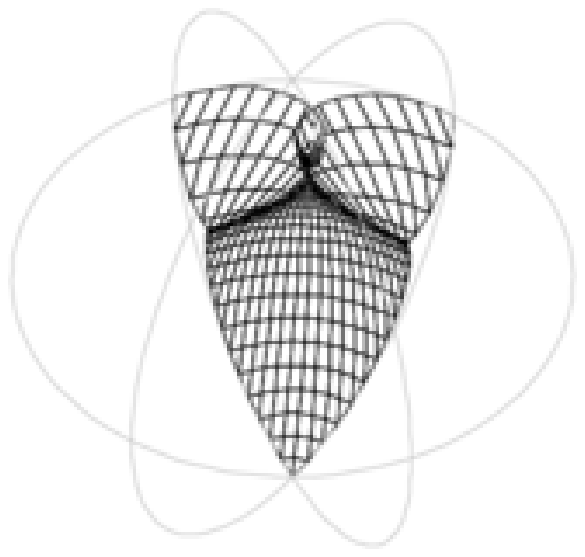} & 
 \includegraphics[width=3.6cm]{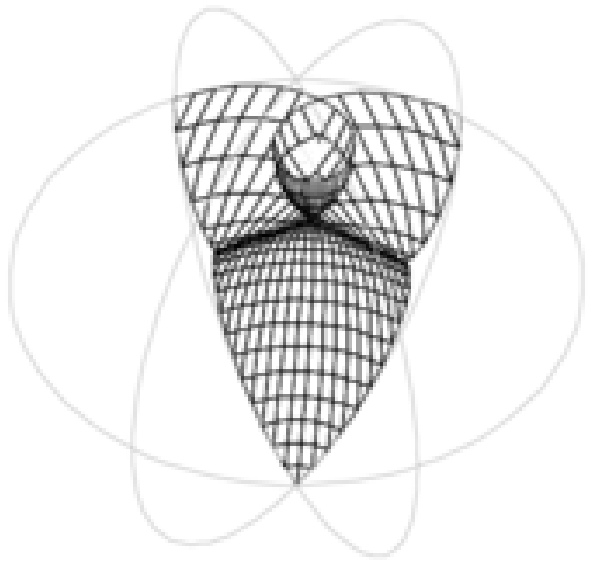} &
 \includegraphics[width=3.6cm]{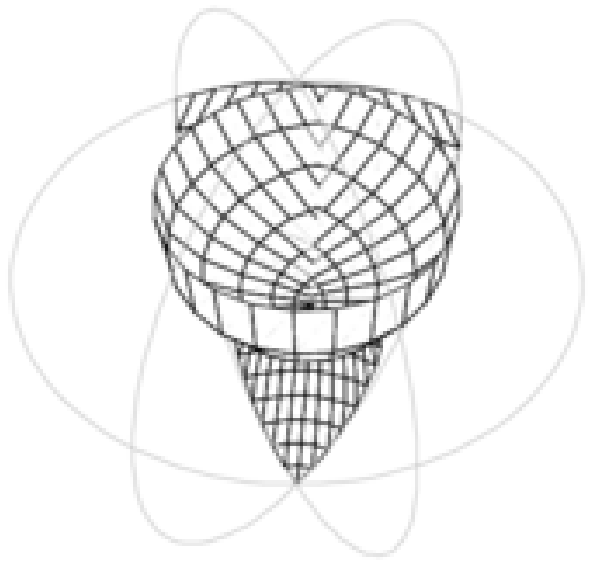} 
\end{tabular}
\end{center}
\caption{Parallel family of the images when $(a,b,c)=(1/6,-1/6,1/2)$}
\label{dihed2}
\end{figure}

\subsection{View of the parallel family}
Relying on the previous discussion, we draw pictures in the case
where $(a,b,c)=(1/6, -1/6, 1/2)$, 
i.e. the case where the monodromy
group is the dihedral group of order 12; refer to \cite{SYY}.

First, the domain $S(X_+)$ is a fan with center at $(0,0)$ of radius $1$
with angle from $0$ to $\pi/3$.
Its image in the sphere, $\partial\BB^3$, is drawn in the left
of Figure \ref{dihed1}. 
The right two figures are the image $HS(X_+)$ and $\chi(DS(X_+))$; 
here, one half of $\chi(DS(X_+)$
is drawn for the sake of a better view.

Second, the parallel family to $HS(X_+)$ is drawn
in Figure \ref{dihed2}: from left to right
and from top to down, the figures gradually change their shape from
$\chi(S(X_+))$ to $\chi(DS(X_+))$.

\begin{figure}
\begin{center}
\includegraphics[width=5cm]{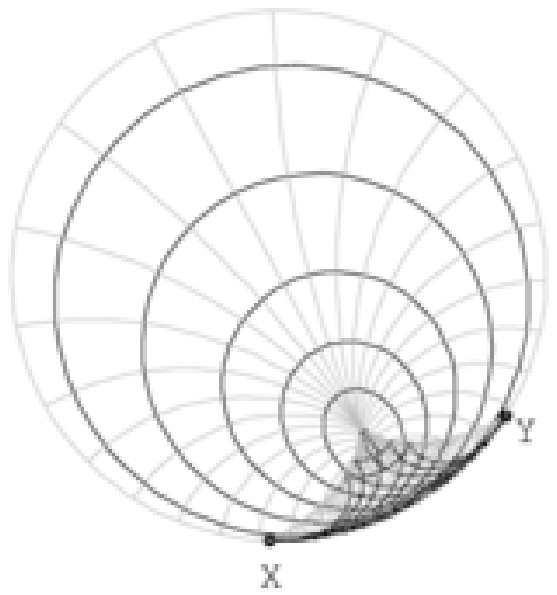}
\hskip1cm 
\includegraphics[width=4.5cm]{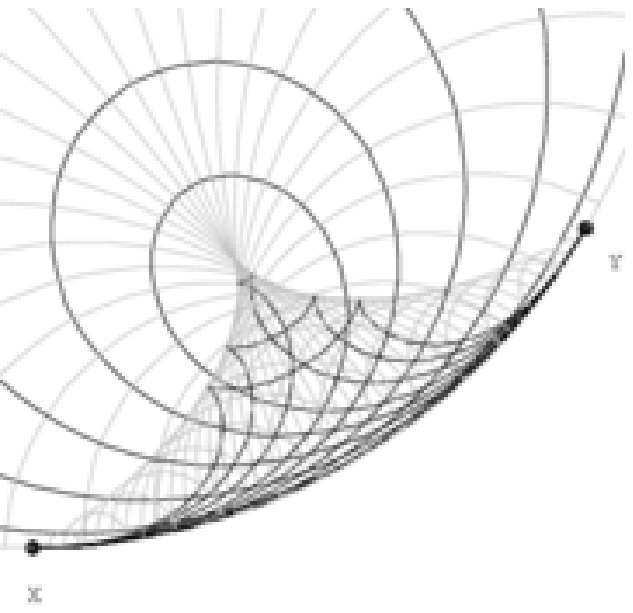} \\
\end{center}
\caption{Left:The section by the equatorial plane of the parallel family,
Right:Caustic locus}
\label{dihed3}
\end{figure}

\begin{figure}[h]
\begin{center}
\includegraphics[width=5cm]{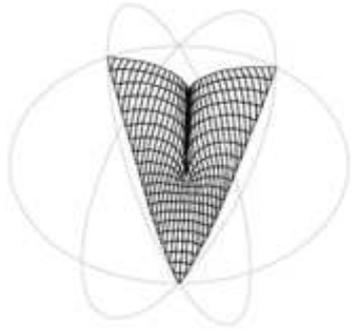} \\
\end{center}
\caption{View of the caustic locus}
\label{dihed4}
\end{figure}

Third, Figure \ref{dihed3} (left) draws the section by the equatorial
plane of the parallel flat fronts.
The gray curves are geodesics joining 
$\chi(S(X_+))$ and $\chi(DS(X_+))$. 
Here, $\chi(S(X_+))$ is the arc $\stackrel{\frown}{XY}$, 
and $\chi(DS(X_+))$ is the arc starting at $X$, going around the circle
once and then terminating at $Y$; 
thus the arc $\stackrel{\frown}{XY}$ is covered twice.
Figure \ref{dihed3}(Right) is an enlargement of the left. 
The enveloping curve of geodesics is the section of the caustic locus,
where some of sections of the flat fronts have cuspidal edges.
Figure \ref{dihed4} gives a total view of the caustic locus.


\end{document}